\documentclass[12pt, leqno]{article}
\usepackage{amsmath}
\usepackage{amssymb}
\usepackage{amscd}
\usepackage{amsthm}
\def\overset#1#2{{\mathrel{\mathop {{#2}_{}}\limits^{#1}}}}
\def\underset#1#2{{\mathrel{\mathop {{}_{} {#2}}\limits_{{#1}_{}}}}}
\def\upplim_#1{\underset{#1}{\overline\lim}\;}
\def\lowlim_#1{\underset{#1}{\underline\lim}\;}

\newcommand{\C}{{\mathbf{C}}}
\newcommand{\cI}{{\mathcal{I}}}

\newcommand{\cS}{{\mathcal{S}}}

\newcommand{\del}{{\partial}}
\newcommand{\delbar}{{\bar\partial}}
\newcommand{\D}{{\mathrm{D}}}

\newcommand{\ex}{{\empty^\exists}}

\newcommand{\fa}{\empty^\forall}

\newcommand{\gm}{\frak{m}}

\newcommand{\iso}{\cong}

\newcommand{\N}{\mathbf{N}}
\newcommand{\nb}{{\bf N.B.}  }

\renewcommand{\O}{{\mathcal{O}}}

\newcommand{\PD}{{{\mathrm{P}\Delta}}}

\newcommand{\R}{{\mathbf{R}}}

\newcommand{\sF}{\mathcal{F}}

\newcommand{\sI}{\mathcal{I}}

\newcommand{\sU}{\mathcal{U}}

\newcommand{\tensor}{{\otimes}}

\newtheorem{lem}[equation]{\bf Lemma}
\newtheorem{thm}[equation]{\bf Theorem}
\newtheorem{claim}[equation]{\it \quad Claim}

\newtheorem{defin}[equation]{\indent {\it Definition}\rm}

\numberwithin{equation}{section}
\newenvironment{pf}
{{\it Proof}\hskip10pt} {\hfill{\it q.e.d.}\par\vskip+10pt}
\title{Another Direct Proof of Oka's Theorem (Oka IX)}
\author{Junjiro {\sc Noguchi}}
\date{Version 1, 2011/aug/10}
\begin{document}
\parskip+2pt
\maketitle
\thispagestyle{empty}
\begin{abstract}
In 1953 K. Oka IX solved in first and in a final form
Levi's problem (Hartogs' inverse problem) for domains or
Riemann domains over $\C^n$ of arbitrary dimension.
Later on a number of the proofs were given;
cf.\ e.g., Docquier-Grauert's paper in 1960, R. Narasimhan's paper
in 1961/62, Gunning-Rossi's book, and H\"ormander's book.
%(in which the holomorphic separability is pre-assumed in the definition
%of Riemann domains and thus the assumption is stronger than the one
%in the present paper).
Here we will give another direct elementary proof of Oka's Theorem,
relying only on Grauert's finiteness theorem by
the {\it induction on the dimension} and the {\it jets over
Riemann domains};  hopefully, the proof is most comprehensive.
\end{abstract}

\section{\large Introduction.}
In 1953 K. Oka \cite{oka} IX solved in first and in a final form
Levi's problem (Hartogs' inverse problem) for domains or
Riemann domains over $\C^n$ of arbitrary dimension
(cf.\ below for notation):

\begin{thm}{\rm (Oka \cite{oka} IX, ('43)/'53\footnote{
It is now possible to confirm that Oka IX published in 1953
was written in French from his notes in Japanese
dated 1943. Cf.\ the introduction of Oka IX, and also
Oka VI published in 1942; see
http://www.lib.nara-wu.ac.jp/oka/index\_eng.html.})}
\label{oka}
Let $\pi:X \to \C^n$ be a Riemann domain, and let
$\delta_{\PD}(x, \del X)$ denote the boundary distance function
with respect to a polydisc $\PD$.
If $- \log \delta_{\PD}(x, \del X)$ is plurisubharmonic,
then $X$ is Stein.
\end{thm}

Besides Oka's original proof there are known
a number of the proofs in generalized forms; e.g.,
Docquier-Grauert \cite{dg}, Narasimhan \cite{nar},
Gunning-Rossi \cite{GR}, and H\"ormander \cite{hoe}
(in which the holomorphic separability is pre-assumed in the definition
of Riemann domains and thus the assumption is stronger than the one
in the present paper).

Here we will give another direct elementary proof of
Oka's Theorem \ref{oka} by making use of the followings
in an essential way, and it is new in this sense
(see the proof of Lemma \ref{key}).
\begin{enumerate}
\itemsep-2pt
\item
The induction on the dimension $n=\dim X$.
\item
The jets over $X$.
\item
Grauert's Finiteness Theorem \ref{gthm} over a strongly
pseudoconvex domain $\Omega$ of a complex manifold
applied not only for the structure sheaf $\O_\Omega$,
but also for a coherent ideal sheaf $\sI \subset \O_\Omega$
(cf.\ Narasimhan \cite{nar}, Docquier-Grauert \cite{dg},
Gunning-Rossi \cite{GR}).
\end{enumerate}
The others are the vanishing of higher cohomologies of coherent
sheaves on polydiscs and on Stein manifolds,
and a sort of $\epsilon$-$\delta$ arguments, to say,
a content presented in Chap.\ 2 of H\"ormander \cite{hoe}
(see, e.g. the proof of Lemma \ref{SR}).
Thus, the proof is elementary, self-contained and
hopefully most comprehensive.

To be precise we give the exact definitions of notions we will use.
\begin{defin}\rm
{\rm (Stein manifold)}
\label{stein}
A connected complex manifold $M$ with the second countability axiom
is called a {\it Stein manifold} if it satisfies the following three
conditions.
Here, $\O(M)$ denotes the set of all holomorphic functions on $M$.
\begin{enumerate}
\item
(Holomorphic separability) For distinct two points $x, y \in M$
there exists an element$f \in \O(M)$ such that $f(x) \not= f(y)$.
\item
(Holomorphic local coordinates)
For an arbitrary point $x \in M$ there are $n$ $(=\dim M)$
elements $f_j \in \O(M)$、 $1 \leq j \leq n$ such that
$(f_j)_{1 \leq j \leq n}$ gives rise to a holomorphic local coordinate
system in a neighborhood of $x$.
\item
(Holomorphic convexity) For a compact subset $K \Subset M$
its holomorphic convex hull
$$
\hat{K}_M=\{x \in M; |f(x)| \leq \max_K |f|, \: \fa{f} \in \O(M)\}
$$
is also compact in $M$.
\end{enumerate}
\end{defin}

\nb In a number of references the definition of Stein manifolds
consists of the above (iii) and the following
{\it K-completeness} due to Grauert \cite{gr55}:
\begin{enumerate}
\item[{(K)}]
``For every point $x \in M$ there exist finitely many
$f_j \in \O(M)$, $1 \leq j \leq l$ such that
all $f_j(x)=0$ and $x$ is isolated in the analytic subset
$\{f_j=0; 1 \leq j \leq l\}$.''
\end{enumerate}
\noindent
In fact, they are equivalent: it is trivial that the present
definition~\ref{stein}
implies the above (K), but the converse is {\em not} trivial
at all (cf.\ Grauert [55]).
\medskip

Let $X$ be a complex manifold and let $\pi:X \to \C^n$ be
a holomorphic map.
\begin{defin}\rm (Riemann domain)
$\pi:X \to \C^n$ or simply $X$ is called a {\it Riemann domain}
if the following properties are satisfied:
\begin{enumerate}
\item
$X$ is connected.
\item
For every point $x \in X$ there are neighborhoods
$U \ni x$ in $X$ and $V \ni \pi(x)$ in $\C^n$ such that
the restriction $\pi|_U: U \to V$ is biholomorphic.
\end{enumerate}
\end{defin}

{\nb} (i) A Riemann domain $X$ is metrizable and hence
$X$ satisfies the second countability axiom.

(ii) In the above definition we do {\em not} assume the holomorphic
separability for a Riemann domain.

A Riemann domain $\hat{\pi}: \hat{X} \to \C^n$ is called a
{\it holomorphic extension} of a Riemann domain $\pi: X \to \C^n$
if there is a holomorphic injection $\iota: X \to \hat{X}$
satisfying
\begin{enumerate}
\itemsep-2pt
\item
$\pi=\hat{\pi} \circ \iota$;
\item
every holomorphic function $f \in \O(X)$ is analytically
continued to an element $\hat{f} \in \O(\hat{X})$.
\end{enumerate}

A Riemann domain $X$ is called a {\it domain of holomorphy} if
there exists no holomorphic extension of $X$ other than $X$ itself.

In this paper $X$ denotes always a Riemann domain.
We take a polydisc $\PD=\PD(0; r_0)$
($r_0=(r_{0j})$) with center at the origin $0 \in \C^n$.
Then by definition there are $\rho>0$ and a neighborhood
$U_\rho(x) \ni x$ for every $x\in X$ such that
$$
\pi|_{U_\rho(x)}: U_\rho(x) \to \pi(x)+\rho\PD
$$
is biholomorphic. The supremum of such $\rho>0$
$$
\delta_{\PD}(x, \del X)=\sup\{\rho>0; \ex U_\rho(x)\} \leq \infty
$$
is called the {\it boundary distance function} of $X$
to the relative boundary.

If $\delta_{\PD}(x, \del X)=\infty$, then $\pi$ is a
holomorphic isomorphism, and thus there is nothing to discuss more.
Henceforth we assume
$\delta_{\PD}(x, \del X)<\infty$ in what follows.

For a subdomain $\Omega \subset X$ we define similarly
$$
\delta_{\PD}(x, \del \Omega)=\sup\{\rho>0; \ex U_\rho(x) \subset \Omega\}.
$$
The boundary distance functions
$\delta_{\PD}(x, \del X)$ and $\delta_{\PD}(x, \del \Omega)$
are continuous with Lipschitz' condition.
For a subset set $A \subset X$ (resp. $A \subset \Omega$) we set
\begin{align*}
\delta_{\PD}(A, \del X)&=\inf_{x \in A} \delta_{\PD}(x, \del X)\\
(\hbox{resp. } \delta_{\PD}(A, \del \Omega) 
&=\inf_{x \in A} \delta_{\PD}(x, \del
 \Omega)).
\end{align*}

{\it Acknowledgment.}  During the preparation of this paper
the author had a number of discussions on K. Oka's works
with Professors K. Kazama, H. Yamaguchi, and S. Hamano,
which were very helpful and of pleasure.
The author would like to express sincere gratitude to
all of them.

\section{\large Preliminaries.}
Here we list up the lemmas and theorems we will use.

\begin{lem}
\label{2.1}
Let $\pi: X \to \C^n$ be a domain of holomorphy, let
$K \Subset X$ be a compact subset, and let
$f \in \O(X)$. If
$$
\delta_{\PD}(x, \del X) \geq |f(x)|, \quad x \in K,
$$
then
$$
\delta_{\PD}(x, \del X) \geq |f(x)|, \quad x \in \hat{K}_X .
$$
In particular, taking $f$ to be constant we have
\begin{equation}
\label{bdrydist}
\delta_{\PD}(K, \del X) =
\delta_{\PD}(\hat{K}_X, \del X).
\end{equation}
\end{lem}

The proof is the same as in the case of univalent domains.
This lemma implies the following as well:
\begin{thm}
If $X$ is a domain of holomorphy, then
$- \log \delta_{\PD}(x, \del X)$ %called the logarithmic boundary function
is plurisubharmonic.
\end{thm}

In general, a complex manifold $M$ is said to be {\it pseudoconvex} if
$M$ carries a plurisubharmonic exhaustion function.
The following is not trivial but elementary due to Oka \cite{oka} IX
(cf.\ Nishino \cite{ni}, p.\ 350):
\begin{lem}
\label{psdc}
If $- \log \delta_{\PD}(x, \del X)$ is plurisubharmonic
(for one fixed $\PD$), then $X$ is pseudoconvex.
\end{lem}

\begin{thm}{\rm (Oka's Fundamental Theorem, I, II, VII, VIII)}
\label{coh}
Let $\PD(0; r)$ be an arbitrary polydisc, and let
$\cI \subset \O_\Omega^N$ be a coherent sheaf of submodules.
Then
$$
H^q(\PD(0; r), \cI)=0 , \qquad q \geq 1.
$$
\end{thm}

This theorem over polydiscs together with Oka's Jok\^uiko\footnote{
A direct English translation may be ``transformation to the upper
space''.
It is a method to imbed the domain under consideration into
a higher dimensional polydisc $\PD$, to extend the analytic objects
over $\PD$, and to solve the problem over $\PD$ by the simplicity of
the space $\PD$. This method was developed by
K. Oka \cite{oka} I$\sim$III and was a very key to
solve Cousin Problems I and II.
}
leads to the following:
\begin{thm}{\rm (Oka-Cartan)}
\label{oc}
Let $M$ be a Stein manifold, and let
$\cS \to M$ be a coherent sheaf.  Then
$$
H^q(M, \cS)=0, \qquad q \geq 1.
$$
\end{thm}

\begin{lem}
\label{2.5}
\begin{enumerate}
\item
Let $\Omega_1 \Subset \Omega_2 \Subset \Omega_3 \Subset X$
be a series of subdomains.
Assume that $\Omega_3$ is Stein.
If
$$
\delta_{\PD}(\del \Omega_1, \del \Omega_3)>
\max_{x \in \del \Omega_2} \delta_{\PD}(x, \del \Omega_3),
$$
then there is an $\O(\Omega_3)$-analytic polyhedron $P$
such that
$$
\Omega_1 \Subset P \Subset \Omega_2 .
$$
\item
An arbitrary holomorphic function $f \in \O(P)$
can be approximated uniformly on compact subsets by
elements of $\O(\Omega_3)$; that is,
$(P, \Omega_3)$ is a Runge pair.
\end{enumerate}
\end{lem}

\begin{pf}
(i) The assumption and \eqref{bdrydist} imply that
$\widehat{(\bar{\Omega}_1)}_{\Omega_3} \Subset \Omega_2$,
and hence such $P$ exists.

(ii) By Theorem \ref{coh} we can apply Oka's Jok\^uiko
to reduce the domain to a polydisc, and is proved.
\end{pf}

Let $\Omega \Subset M$ be a relatively compact domain.
\begin{defin}\rm
$\Omega$ is said to be {\it strongly pseudoconvex}
if there are a neighborhood $U (\subset M)$ of the boundary
$\del \Omega$ of $\Omega$, and a real valued $C^2$ function
$\phi: U \to \R$ satisfying the conditions\maketitle
\begin{enumerate}
\item
$\{x \in U : \phi(x)<0\}=\Omega \cap U$,
\item
$i \del \delbar \phi(x)>0$ ($x \in U$).
\end{enumerate}
\end{defin}

\begin{thm}{\rm (Grauert \cite{gr58}, \cite{gr62})}
\label{gthm}
Let $\Omega \Subset M$ be a strongly pseudoconvex domain.
Let $\sF$ be a coherent sheaf defined over a neighborhood
of the closure $\bar{\Omega}$. Then we have
$$
\dim H^q(\Omega, \sF)< \infty, \quad q \ge 1.
$$
\end{thm}

We will use this theorem for the structure sheaf and an
ideal sheaf of a closed complex submanifold.
In the first, we apply this for
$\sF=\O_M$ to deduce
\begin{thm}
\label{gthm1}\rm
Let $\Omega$ be as in Theorem \ref{gthm}.
Then $\Omega$ is holomorphically convex.
\end{thm}

{\nb} The above described was the circumstance just after
Grauert \cite{gr58} ('58), and before 
Docquier-Grauert \cite{dg} ('60) and Narasimhan \cite{nar} ('61/'62).

\section{\large A Proof of Oka's Theorem \ref{oka}.}

By Lemma \ref{psdc} it suffices to show the following for the proof.
\begin{thm}
A pseudoconvex Riemann domain is Stein.
\end{thm}

Under the assumption we take a plurisubharmonic exhaustion function
$\phi:X \to [-\infty, \infty)$.
The following lemma is our key.
\begin{lem}
\label{key}
If $\Omega \Subset X$ is a strongly pseudoconvex domain,
then $\Omega$ is Stein.
\end{lem}
\begin{pf}
We use the induction on the dimension $n\geq 1$.

{\bf (a) } $n=1$: In this case ${\Omega}$ is an open Riemann surface
and hence by Behnke-Stein's Theorem it is Stein.
For the completeness we show this with the preparation in \S2.
The holomorphic convexity is finished by Theorem \ref{gthm1}.
The holomorphic local coordinates follow just from the definition
of Riemann domain. It is remaining to show the holomorphic separability.

Take two distinct points $a, b \in \Omega$.
If $\pi(a)\not=\pi(b)$, the proof is done.
Suppose that $\pi(a)=\pi(b)$. By a translation of $\C$
we may assume that $\pi(a)=\pi(b)=0 \in \C$.
Let $U_0 \ni a$ be a neighborhood such that $U_0 \not\ni b$
and $\pi|_{U_0}: U_0 \to \Delta(0; \delta)$ with $\delta>0$
is biholomorphic.
Put $U_1=\Omega \setminus \{a\}$.
Then $\sU=\{U_0, U_1\}$ is an open covering of $\Omega$.
For each $k \in \N$ we set
$$
\gamma_k(x)=\frac{1}{\pi(x)^k},\quad x \in U_0 \cap U_1.
$$
Then $\gamma_k$ defines an element of $H^1(\sU, \O_\Omega)$.
It is noted that
$H^1(\sU, \O_\Omega) \hookrightarrow H^1(\Omega, \O_\Omega)$
is injective.
By Theorem \ref{gthm} there is a non-trivial linear relation
$$
\sum_{k=1}^{h} c_k \gamma_k=0, \quad c_k \in \C, \, c_h\not=0 .
$$
Therefore there are elements $f_j \in \O(U_j), j=0,1$
such that
$$
f_1(x) - f_0(x)= \sum_{k=1}^{h} c_k \frac{1}{\pi(x)^k},
\quad x \in U_0 \cap U_1.
$$
Thus we obtain a meromorphic function
in $\Omega$ with a pole only at $a$,
$$
F=f_1=f_0 + \sum_{k=1}^{h} c_k \frac{1}{\pi^k} .
$$
From the construction we get
\begin{align*}
\pi(x)^h F(x) & \in \O(\Omega),\\
\pi(a)^h F(a) & =c_h \not= 0,\\
\pi(b)^h F(b) & =0.
\end{align*}
Therefore $a$ and $b$ are separated by an element of $\O(\Omega)$.

{\bf (b) } We assume the assertion holds in $\dim X=n-1$.
Let $\dim X=n\geq 2$.
By the definition of Riemann domain it is sufficient to prove
the holomorphic convexity and the holomorphic separability;
the first is finished by Theorem \ref{gthm1}, and the latter
remains to be shown.

{\bf (1) } We take arbitrary distinct points $a,b \in \Omega$.
As in (a) we may assume that $\pi(a)=\pi(b)=0$.
Taking a hyperplane $L=\{z_n=0\}$, we consider the restriction
$$
\pi_{X'}: X'=\pi^{-1}L  \to L.
$$
Since $L \iso \C^{n-1}$ (biholomorphic),
every connected component $X''$ of $X'$ is
$(n-1)$ dimensional Riemann domain.
The restriction $\phi|_{X''}$ is a plurisubharmonic exhaustion
function. By the induction hypothesis $X''$ is Stein.

{\bf (2) }  Let $\gm \langle a \rangle \subset \O_{X', a}$ be the
maximal ideal of the local ring $\O_{X', a}$ and let
$\gm_a^k$ denote the $k$-th power.  Set
$$
\gm^k \langle a, b \rangle=
\gm^k \langle a \rangle \tensor \gm^k \langle b \rangle
\subset \O_{X'}.
$$
This is a coherent ideal sheaf of $\O_{X'}$.

Since every connected component of $X'$ is Stein,
Theorem \ref{oc} implies the existence
of $g_k \in \O(X')$ for each $k \in \N$ such that
\begin{align}
\label{jet}
\underline{g_k}_a &\equiv 0 \quad (\hbox{mod } \gm^{k-1} \langle a, b
 \rangle_a),\\
\nonumber
\underline{g_k}_a & \not\equiv 0 \quad (\hbox{mod } \gm^{k} \langle a, b
 \rangle_a),\\
\nonumber
\underline{g_k}_b &\equiv 0 \quad (\hbox{mod } \gm^{k} \langle a, b
 \rangle_b),
\end{align}
where $\underline{g_k}_a$ stands for a germ of $g_k$ at $a$.

{\bf (3) }  We put $\Omega'=\Omega \cap X'$.
Let $\cI$ be the ideal sheaf of the analytic subset
$X' \subset X$.
By Oka's Second Coherence Theorem (\cite{oka} VII, VIII) $\cI$
is coherent.\footnote{\, There seems to be a confusion in the historical
comprehension of the developement of the ``coherence theorems''.
 In Oka VII and VIII K. Oka proved
 {\it three fundamental coherence theorems}.
Firstly in Oka VII which was received in 1948 and published
in 1950, he proved the coherence of the structure sheaf $\O_{\C^n}$
on $\C^n$ (Oka's First Coherence Theorem),
and he was writing in two places that in the
forthcoming paper he would deal with the coherence of
ideal sheaves of analytic subsets,
 ``{\it id\'eaux g\'eom\'etriques de domaines ind\'etermin\'es}''
he termed,
and that one would see it to hold without any assumption;
see 1) the last six lines of the paper at p.\ 27,
and 2) the last two lines of p.\ 7 to the line just before \S3 of p.\ 8.
There he wrote that there are two cases for which the coherence problem
are solvable, the first is that of $\O_{\C^n}$ dealt with in VII,
%(Oka's First Coherence Theorem),
and the second is that of the ideal sheaf
of an analytic subsets (Oka's Second Coherence Theorem),
of which proof appeared in Oka VIII in 1951, while H. Cartan's proof
appeared in 1950 in the same volume as Oka VII,
to which the theorem is attributed in most references.

For this many refer only to the first point 1), but never to the second
point 2) so far by the knowledge of the present author,
where K. Oka was writing more detailed descriptions what should be done
for the second coherence theorem.
In VIII he wrote its proof % of the second coherence 
and moreover
proved the coherence of normalizations (Oka's Third Coherence Theorem).
For a convenience we give a complete list of of K. Oka's paper
at the end of the references, which is not very long but hard to find
a correct one.
}
Restriction this to $\Omega$
we have a short exact sequence:
$$
0 \to \cI \to \O_\Omega \to \O_{\Omega'} \to 0.
$$
This implies the following exact sequence,
\begin{equation}
\label{ex}
\O(\Omega) \to \O(\Omega') \ \overset{\delta}{\to} \
 H^1(\Omega, \cI).
\end{equation}
We write $g_k$ for the restriction of $g_k$ to $\Omega'$ by the same
 letter.
We have that $\{\delta(g_k)\}_{k \in \N} \subset H^1(\Omega, \cI)$.
By Theorem \ref{gthm} $H^1(\Omega, \cI)$ is finite dimensional,
and thus there is a non-trivial linear relation
$$
\sum_{k=k_0}^N c_k \delta(g_k)=0, \quad c_k \in \C, \ N< \infty .
$$
We may assume that $c_{k_0}\not=0$. It follows from \eqref{ex}
that there is an element $f \in \O(\Omega)$ such that
$$
f|_{\Omega'}=\sum_{k=k_0}^N c_k g_k.
$$
We use $\pi=(z_1, \ldots, z_n)$ as a holomorphic local coordinate
system in a sufficiently small neighborhood of $a \in \Omega$,
$z'=(z_1, \ldots, z_{n-1})$. Then we get
\begin{equation}
\label{sep1}
f(z)= \sum_{k=k_0}^N c_k g_k(z') + h(z)\cdot z_n,
\end{equation}
where $h(z)$ is a holomorphic function in a neighborhood of $a$.
It follows from \eqref{jet} that there is a partial
differentiation of order $k_0$ in $z'$
$$
{\D}=\frac{\del^{k_0}\qquad\qquad }{\del z_1^{\alpha_1} \cdots \del
 z_{n-1}^{\alpha_{n-1}}},
\quad \sum_{j=1}^{n-1} \alpha_j=k_0
$$
such that
\begin{align}
\label{sep2}
\D g_{k_0}(a)&\not=0,\\
\nonumber
\D g_{k}(a) &=0, \quad k >k_0,\\
\nonumber
\D g_k(b) &= 0, \quad k \geq k_0.
\end{align}
The definition of $\D$ and \eqref{sep1} imply that
$$
\D f(z)= \sum_{k=k_0}^N c_k \D g_k(z') + (\D h(z)) \cdot z_n.
$$
Since $z_n=0$ at $a$ and $b$, \eqref{sep2} leads to
$$
\D f(a) \not= 0, \quad \D f(b)=0.
$$
Since $\D f \in \O(\Omega)$, the holomorphic separability of $\Omega$
was proved.
\end{pf}

We set
$$
X_c= \{x \in X; \phi(x)<c\}, \quad c \in \R .
$$
For $X$ being Stein it suffices to prove the followings:
\begin{lem}
\label{SR}
\begin{enumerate}
\item
$X_c$ is Stein for an arbitrary $c \in \R$;
\item
For every pair of $c < b$, $(X_c, X_b)$ is a Runge pair.
\end{enumerate}
\end{lem}
\begin{pf}
(i) Let $K \Subset X_c$ be a compact subset.  We put
$$
\eta=\delta_{\PD}(K, \del X_c)\, (>0) .
$$
We take $b>c$ so that
\begin{equation}
\label{3.8}
\max_{x \in \del X_c} \delta_{\PD}(x, \del X_b)< \eta .
\end{equation}
Since $\|\pi(x)\|^2$ is strongly plurisubharmonic everywhere
and $\phi$ is plurisubharmonic,
there exists a strongly pseudoconvex domain $\Omega$ such that
$$
X_c \Subset \Omega \Subset X_b .
$$
By Lemma \ref{key} $\Omega$ is Stein.
Therefore conditions (i) and (ii) of Definition \ref{stein}
are satisfied, and there remains (iii) (holomorphic convexity)
to be shown.
\begin{claim}
\label{3.9}
$\hat{K}_{X_c} \Subset X_c$ .
\end{claim}
$\because \,)$ The application of \eqref{bdrydist}
to $K \Subset \Omega$ yields
$$
\delta_{\PD}(\hat{K}_\Omega, \del \Omega)
=\delta_{\PD}(K, \del \Omega) > \eta.
$$
On the other hand, from \eqref{3.8} it follows that
$$
\max_{x \in \del X_c} \delta_{\PD}(x, \del \Omega)< \eta.
$$
The above two equations imply
\begin{equation}
\label{3.10}
\hat{K}_{X_c} \subset \hat{K}_\Omega \Subset X_c \, .
\end{equation}

(ii) We use the same notation as in (i).

{\bf (1) }  We now know that all $X_c$ $(c \in \R)$ are Stein.
Therefore, replacing $\Omega$ by $X_b$ in the above
arguments in (i), we see that
\begin{equation}
\label{3.11}
\hat{K}_{X_c} \subset \hat{K}_{X_b} \Subset X_c \Subset X_b \, .
\end{equation}

\begin{claim}
\label{3.12}
  \quad $\hat{K}_{X_c} = \hat{K}_{X_b}$.
\end{claim}
$\because \,)$  By \eqref{3.11} we can take an $\O(X_b)$-analytic
polyhedron $P$ such that
$$
\hat{K}_{X_c} \subset \hat{K}_{X_b}
\Subset P \Subset X_c \Subset X_b \, .
$$
If there is a point $\zeta \in \hat{K}_{X_b} \setminus \hat{K}_{X_c}$,
then there is some $g \in \O(X_c)$ such that
$$
\max_K |g| < |g(\zeta)|.
$$
By Lemma \ref{2.5} (ii) $g$ can be approximated uniformly on
$\hat{K}_{X_b}$ by an element of $\O(X_b)$.
Hence there is a holomorphic function $f \in \O(X_b)$ such that
$$
\max_K |f| < |f(\zeta)|.
$$
This is absurd.

{\bf (2) }  It follows from Claim \ref{3.12} that
\begin{equation}
\label{3.13}
\hat{K}_{X_c} = \hat{K}_{X_t}, \quad c \leq \ \fa{t} \leq b.
\end{equation}
We set
$$
E= \{t \geq c\, ;  \hat{K}_{X_t}=\hat{K}_{X_c}\} \subset [c, \infty) .
$$
By definition $t \in E$ implies $[c, t] \subset E$.
The result of (1) shows that $E$ is an open subset of $[c, \infty)$.

{\bf (3) }  We put 
$a=\sup \, E$.

\begin{claim}
\quad $a=\infty$; i.e., $E=[c, \infty)$.
\end{claim}
$\because \,)$  Suppose that $a< \infty$.
From the definition we obtain
$$
K_1=\hat{K}_{X_c} = \hat{K}_{X_t}, \quad c \leq \ \fa{t} <a .
$$
Letting $t<a$ sufficiently close to $a$, we have
$$
\delta_{\PD}(K_1, \del X_a)> \max_{x \in \del X_t}
\delta_{\PD}(x, \del X_a).
$$
Because $X_a$ is Stein,
$$
\delta_{\PD}(\hat{K_1}_{X_a}, \del X_a)=
\delta_{\PD}(K_1, \del X_a)> \max_{x \in \del X_t}
\delta_{\PD}(x, \del X_a).
$$
Thus, $\hat{K_1}_{X_a} \Subset X_t$ follows.
One gets
$$
\hat{K}_{X_t} \subset \hat{K}_{X_a}
\subset \hat{K_1}_{X_a} \Subset X_t \Subset X_a.
$$
In the same way as in (1) we see that $\hat{K}_{X_t} = \hat{K}_{X_a}$.
Therefore, $a \in E$.
Since $E$ is open, there exists a number $a' \in E$ with
$ a'>a$. This contradicts to the choice of $a$.

{\bf (4) } It follows from (2) that for arbitrary
$c<b$ and a compact subset $K \Subset X_c$,
$$
\hat{K}_{X_c}=\hat{K}_{X_b}.
$$
Therefore, Oka's Jok\^uiko and Theorem \ref{coh} imply that
$(X_c, X_b)$ is a Runge pair.
\end{pf}

\bigskip

\rightline{Graduate School of Mathematical Sciences}
\rightline{The University of Tokyo}
\rightline{Komaba, Meguro,Tokyo 153-8914}
\rightline{e-mail: noguchi@ms.u-tokyo.ac.jp}
\end{document}